\documentclass{amsart}
\usepackage{amssymb,euscript,amsmath, mathrsfs}
\usepackage[dvips]{graphicx}
\usepackage[dvips]{color}

\newcounter{ENUM}
\newcommand{\itm}{\item}
\newenvironment{ilist}{\renewcommand{\theENUM}{\roman{ENUM}}\renewcommand{\itm}{\addtocounter{ENUM}{1}\item[(\theENUM)]}\begin{itemize}\setcounter{ENUM}{0}}{\end{itemize}}
\newenvironment{Ilist}{\renewcommand{\theENUM}{\Roman{ENUM}}\renewcommand{\itm}{\addtocounter{ENUM}{1}\item[(\theENUM)]}\begin{itemize}\setcounter{ENUM}{0}}{\end{itemize}}
\newenvironment{alist}[1][0]{\renewcommand{\theENUM}{\alph{ENUM}}\renewcommand{\itm}{\addtocounter{ENUM}{1}\item[\theENUM)]}\begin{itemize}\setcounter{ENUM}{#1}}{\end{itemize}}

\newcommand{\margh}[1]{}

\input xy
\xyoption{all}
\CompileMatrices

\def\PP{{\mathbb P}}

\def\CC{{\mathbb C}}

\def\cC{{\mathcal C}}

\def\cH{{\mathcal H}}

\def\cM{{\mathcal M}}

\def\fg{{\mathfrak g}}

\newtheorem{thm}{Theorem}[section]
\newtheorem{prop}[thm]{Proposition}
\newtheorem{lem}[thm]{Lemma}
\newtheorem{cor}[thm]{Corollary}

\theoremstyle{definition}

\newtheorem{ques}[thm]{Question}

\newtheorem{sit}[thm]{Situation}
\theoremstyle{remark}

\newtheorem{rem}[thm]{Remark}

\numberwithin{equation}{section}
\numberwithin{figure}{section}

\begin{document}
\title{The irreducibility of certain pure-cycle Hurwitz spaces}
\author{Fu Liu and Brian Osserman}
\begin{abstract} We study ``pure-cycle'' Hurwitz spaces, parametrizing 
covers of the projective line having only one ramified point over each 
branch point. We start with the case of genus-$0$ covers, using a 
combination of limit linear series theory and group theory to show that 
these spaces are always irreducible. In the case of four branch points, 
we also compute the associated Hurwitz numbers. Finally, we give a 
conditional result in the higher-genus case, requiring at least $3g$ simply 
branched points. These results have equivalent formulations in group theory, 
and in this setting complement results of Conway-Fried-Parker-V\"olklein.
\end{abstract}

\thanks{The first author was supported by fellowships from the Clay 
Mathematics Institute and MSRI, and the second author was partially 
supported by a fellowship from the
National Science Foundation during the preparation of this paper.}
\maketitle

\section{Introduction}

In this paper, we use a combination of geometric and group-theoretic 
techniques to prove a result with equivalent statements in both fields.
The geometric statement is that certain genus-$0$ Hurwitz spaces (the
``pure-cycle'' cases) are always irreducible, while the group-theoretic 
statement is that the corresponding factorizations into cycles always lie 
in a single pure braid group orbit. ``Pure-cycle'' refers to the hypothesis
that for our covers, there is only a single ramified point over each
branch point. The main significance for us of this condition is that it 
allows us to pass relatively freely between the point of view of branched
covers, where one moves the branch points freely on the base curve, and 
linear series, where one moves the ramification points freely on the 
covering curve. This facilitates induction, as it is easier to stay within
the pure-cycle case from the point of view of linear series.

Our result is close to optimal in the 
sense that if one drops either of the pure-cycle or genus-$0$ hypotheses, 
one quickly runs into cases where the Hurwitz spaces have more than one 
component. However, we do prove a conditional generalization to 
higher-genus pure-cycle Hurwitz spaces having at least $3g$ simply branched
points, depending on a positive answer to a different geometric question 
which is closely related to an old question of Zariski.

Our immediate motivation for studying the pure-cycle situation is its
relation to linear series: specifically, if one wishes to prove statements
on branched covers via linear series arguments, the pure-cycle situation
is the natural context to examine. A good understanding of the classical
situation is therefore important to studying other cases, such as that of
positive characteristic. In particular, our main theorem allows for a much
simpler proof of a stronger result in \cite{os12} than would otherwise be
possible. However, we also remark that a good understanding of the components
of Hurwitz spaces has given rise to a wide range of substantial 
applications: the classical proof of Severi that $\cM_g$ is connected
\cite{fu3}; number-theoretic applications such as Davenport's problem,
and Thompson's genus-0 problems in group theory, both due to Fried 
\cite{fr2}; and the Fried-V\"olklein description of the absolute 
Galois group of certain fields in inverse Galois theory 
\cite{f-v1},\cite{f-v2}. Furthermore, our results in particular should 
provide good test cases for Fried's Main Conjecture of modular tower theory
\cite{fr5}. See also the survey article \cite{fr4}. 

We now state our results more precisely. We will recall/fix our terminology 
in the next section.

The following proposition is well known, although the equivalence of the 
first two and last two conditions depends heavily on the fact that we 
restrict our attention to covers with a single ramified point over each
branch point. We will recall the proof in the following section.

\begin{prop}\label{equivs} Given $d$ and $\vec{e}=(e_1,\dots,e_r)$ with 
$2d-2=\sum_i(e_i-1)$, the following are equivalent:
\begin{alist}
\itm the Hurwitz factorizations for $(d,r,0,\vec{e})$ all lie in a single 
orbit of the pure braid group.
\itm the space $\cH(d,r,0,\vec{e})$ is irreducible, where 
$\cH(d,r,0,\vec{e})$ 
is the Hurwitz space parametrizing $r$ distinct points $Q_1,\dots,Q_r$
on $\PP^1$ together with a genus-0 cover of $\PP^1$, such that each $Q_i$ 
has a single point over it ramified to order $e_i$, and the rest unramified;
\itm the space $MR:=MR(\PP^1,\PP^1, \vec{e})$ is irreducible, where
$MR$ is the space parametrizing $r$ distinct points $P_1,\dots,P_r$ on 
$\PP^1$ together with a rational function $f:\PP^1 \to \PP^1$ of degree $d$ 
and ramified to order $e_i$ at $P_i$ (on the source curve) for all $i$; 
\itm the space $G^1_d:=G^{1}_d(\PP^1,\vec{e})$ is irreducible,
where $G^1_d$ is the space parametrizing $r$ distinct points $P_1,\dots,P_r$
on $\PP^1$ together with a linear series of dimension $1$ and degree $d$, 
having ramification $e_i$ at $P_i$ for all $i$;
\end{alist}
\end{prop}

Our main theorem is then the following:

\begin{thm}\label{main} Given $d,r$ and $e_1,\dots,e_r$ with 
$2d-2=\sum_i(e_i-1)$, the equivalent conditions of Proposition \ref{equivs}
always hold.
\end{thm}

Our proof follows the general structure of Eisenbud and Harris' argument in
\cite{e-h7}, where they prove the irreducibility of certain families of
linear series without prescribed ramification. However, while they work 
exclusively from the perspective of linear series, we have to switch 
back and forth between points of view. 
Starting from the perspective of linear series, we use a degeneration 
argument and the tools of limit linear series to reduce to a base case of 
four points on $\PP^1$, and then solve that case directly, after switching
to the group-theoretic point of view of Hurwitz factorizations. Our explicit
work in the base case also computes the Hurwitz numbers for that case:

\begin{thm}\label{hur-num} Given $d$ and $\vec{e}=(e_1,\dots,e_4)$ with 
$2d-2=\sum_i(e_i-1)$, we have the following formula for the Hurwitz number: 
$$h(d,r,0,\vec{e})=\min\{e_i(d+1-e_i)\}_i.$$  

We can also describe the Hurwitz factorizations in this case completely
explicitly.
\end{thm}

Lastly, in Theorem \ref{high-gen} below, we again use limit linear series
techniques to prove a conditional version of Theorem \ref{main} for 
pure-cycle cases of higher genus having at least $3g$ simply branched 
points, depending on a positive answer to Question \ref{zar-q} below, a 
geometric question closely related to an old question of Zariski.

The higher-genus result could be seen as having the spirit of an effective 
version in the pure-cycle case of results of Conway-Fried-Parker-V\"olklein.
Our main theorem also generalizes a theorem of Fried \cite[Thm.\ 1.2]{fr1}, 
which implies the case of our Theorem \ref{main} in which $e_i=3$ for all 
$i$.

Finally, we remark that the combination of the genus-0 and pure-cycle 
conditions imply that our monodromy groups are always either cyclic, $S_d$, 
or $A_d$; we show this, independently of the proof of our main results, in 
Theorem \ref{nielsen} below.

\section*{Acknowledgements}

We would like to thank Michael Fried, Kay Magaard and David Harbater for 
help with context and references, and Robert Guralnick for assistance with 
the proof of Theorem \ref{nielsen}.

\section{Notation and terminology}

We quickly recall terminology and fix our notation. For geometric 
statements, we assume throughout that we are working over $\CC$. 

Our notation for permutations will always be to express them as products
of cycles. Given $\sigma \in S_d$, we will say that a number 
$k \in \{1,\dots,d\}$ is {\bf in the support} of $\sigma$ if 
$\sigma(k) \neq k$.

Given a permutation $\sigma$ (or conjugacy class $T$) of $S_d$, we define
its {\bf index} $\iota(\sigma)$ as follows:
if $a_{1} \leq a_{2} \leq \dots \leq a_{m}$ is the 
corresponding partition, then $\iota(\sigma):=\sum_{i=1}^m(a_i-1)$. 
We then say that a tuple $(d,r,g,(T_1,\dots,T_r))$ 
constitutes the data of a {\bf Hurwitz problem}, where 
$d\geq 1,r\geq 2,g\geq 0$, the $T_i$ are conjugacy classes in $S_d$, and we
require $2d-2+2g=\sum_i \iota(T_i)$.

Associated to a Hurwitz problem we have the group-theoretic question of
finding all {\bf Hurwitz factorizations} $(\sigma_1,\dots,\sigma_r)$,
where:
\begin{ilist}
\itm $\sigma_i \in T_i$; 
\itm $\sigma_1 \dots \sigma_r=1$;
\itm the $\sigma_i$ generate a transitive subgroup of $S_d$.
\end{ilist}
We say that two Hurwitz factorizations are {\bf equivalent} if they are
related by simultaneous conjugation by an element of $S_d$. 
We call the number of equivalence classes of Hurwitz factorizations 
the {\bf Hurwitz number} $h(d,r,g,(T_1,\dots,T_r))$. See Remark
\ref{rem-nielsen} below for discussion of some relating and conflicting
notation in the literature.

Geometrically, we also have the {\bf Hurwitz space} 
$\cH(d,r,g,(T_1,\dots,T_r))$, 
parametrizing $r$-tuples of marked points on $\PP^1$, together with covers
of degree $d$ and genus $g$, unramified 
away from the marked points, and with monodromy type $T_i$ at the $i$th
marked point for all $i$. For a fixed choice of marked points, such covers
correspond to Hurwitz factorizations up to equivalence, so the degree of
$\cH(d,r,g,(T_1,\dots,T_r))$ over the space $\cM_{0,r}$ parametrizing 
marked points is given by the Hurwitz number.

We say that a Hurwitz problem is {\bf pure-cycle} if each $T_i$ consists of
a single cycle. Throughout this paper, we restrict our attention to 
pure-cycle Hurwitz problems, and we replace the $T_i$ by integers 
$e_i \geq 2$ giving the length of the cycle. We thus have 
$2d-2+2g=\sum_i (e_i-1)$ as the condition on our data.

We recall that the {\bf Artin braid group} $B_r$ acts on tuples 
$(\sigma_1,\dots,\sigma_r)$ in $S_d$ with $\sigma_1 \dots \sigma_r=1$,
preserving the group generated by the $\sigma_i$.
The $i$th generator acts by replacing $(\sigma_i,\sigma_{i+1})$ by
$(\sigma_{i+1},\sigma_{i+1}^{-1} \sigma_i \sigma_{i+1})$. The kernel
of the natural map $B_r \to S_r$ is the {\bf pure braid group}, which not
only preserves $\sigma_1 \dots \sigma_r=1$, but sends each $\sigma_i$ to
a conjugate of itself in the group generated by all the $\sigma_i$. We thus
see that the pure braid group acts on the set of Hurwitz factorizations,
and it is the orbits of this action which we will study.

Note as a consequence of the geometric definition of Hurwitz number that
the number is clearly invariant under reordering of the $e_i$. We can also
see this purely in terms of group theory by making use of the braid group
action to permute the $e_i$ arbitrarily.

We will also be working from the point of view of {\bf linear series},
which from our point of view will always have dimension $1$ and be
basepoint free: in this situation, a linear series of dimension $1$ and
degree $d$ (also called a $\fg^1_d$) on a curve $C$ is simply a map to
$\PP^1$ of degree $d$, considered up to automorphism of the image space.
We remark that the basepoint-free hypothesis will not cause us any problems,
as we will always be working with spaces of linear series with all 
ramification specified.

As a simple case of the sort of analysis we will carry out in the four-point
case, we recall the answer in the case of three points:

\begin{lem}\label{3pts} The Hurwitz number for $(d,3,0,(e_1,e_2,e_3))$
is always $1$, corresponding to the factorization:
\begin{align*}
\sigma_1 &=(d-e_2, d-e_2-1, \dots, 2, 1, e_3, e_3+1, e_3+2, \dots, d-1, d), \\
\sigma_2 &=(d, d-1, \dots, d-e_2+2, d-e_2+1),\text{ and} \\
\sigma_3 &= (1,2,\dots, e_3-1,e_3).
\end{align*}
\end{lem}

\begin{proof} First, note that by transitivity and the fact that 
$\sigma_1 = \sigma_3^{-1} \sigma_2^{-1}$, we have that 
$\sigma_3$ and $\sigma_2$ together act non-trivially on all of 
$\{1,\dots,d\}$, and their actions must therefore overlap on a subset of 
cardinality exactly $e_2+e_3-d=d+1-e_1$.

To complete the proof, one observes
that if a sequence of precisely $m$ consecutive elements in the cycle 
representation of $\sigma_2$ also appear in $\sigma_3$, at most 
$m-1$ of them can remain fixed by $\sigma_2 \sigma_3$.
It follows that in order for $\sigma_2 \sigma_3$ to be an $e_1$-cycle, the 
overlap must form a single contiguous portion of each of $\sigma_2$ and 
$\sigma_3$, from which one easily concludes the desired statement. 
\end{proof}

Finally, we recall:

\begin{proof}[Proof of Proposition \ref{equivs}]
The equivalence of (i) and (ii) is classical and quite
general: the basic idea is that the monodromy cycles of a cover depend
not only on the cover, but also on a choice of local monodromy generators
of the fundamental group of the base; all such choices of generators are
related by braid operations, and each braid operation can be achieved as
monodromy of the Hurwitz space by moving the marked points of the base around
one another. For a slightly different exposition, see 
\cite[Prop.\ 10.14 (a)]{vo1}; note that the situation is slightly different
because he considers Hurwitz spaces with unordered branch points and full
braid orbits, but the argument is the same in our case of ordered branch
points and pure braid orbits.

Similarly, the equivalence of (iii) and (iv) is equally basic: the space
$G^1_d$ is obtained from the space $MR$ simply by modding out by the (free)
action of the automorphism group of the base $\PP^1$, so $MR$ is a 
$PGL_2$-bundle over $G^1_d$, and one space is irreducible if and only if the
other is.

Next, because we have restricted to Hurwitz spaces in which there is a
single ramified point over each branch point, the comparison of $MR$ and
$\cH(d,r,0,\vec{e})$ is almost equally straightforward. First suppose
$r \geq 3$. If we denote by $\widehat{MR}$ the open subscheme of $MR$ for 
which the map $f$ sends the marked ramification points to distinct points, 
then because $r \geq 3$, we have that $\widehat{MR}$ is a $PGL_2$-bundle over 
$\cH(d,r,0,\vec{e})$, so one is irreducible if and only if the other is. 
But then an easy deformation-theory argument shows
that any component of $MR$ dominates the $(\PP^1)^r$
parametrizing the branch points of the map $f$ \cite[Cor.\ 3.2]{os3},
so we see that $\widehat{MR}$ is dense in $MR$, completing the
desired equivalences for irreducibility. Finally, if $r=2$, the only 
maps are, up to automorphism, $x \mapsto x^d$, so it is easy to see that 
both $MR$ and $\cH(d,r,0,\vec{e})$ are irreducible. 
\end{proof}

\begin{rem}\label{rem-nielsen} Our terminology of Hurwitz problem (and more 
specifically, the
associated set of Hurwitz factorizations) is closely related to the
more standard terminology ``Nielsen class'', for which one also specifies a
subgroup $G$ which the $\sigma_i$ must generate, and assigns the $T_i$ as
conjugacy classes within that subgroup. 

The Nielsen class is frequently 
better because it gives a finer combinatorial invariant: the Hurwitz 
factorizations for a given Hurwitz problem are a disjoint union over 
different Nielsen classes, and likewise the Hurwitz space is a disjoint 
union over spaces associated to different Nielsen classes. Our main theorem
immediately implies that for the cases we study, a Hurwitz problem consists
of only a single Neilsen class. See also Theorem \ref{nielsen} below for
a direct proof of this fact. 

One has to be slightly careful in comparing
statements, since the Nielsen class terminology also allows for different 
equivalence relations on the Hurwitz factorizations (for instance, working 
up to inner automorphism of $G$).

We also remark that our terminology of Hurwitz number, although standard in
some areas, conflicts with the usage in \cite{fr3}.
Specifically, in {\it loc.\ cit.}, the term ``Hurwitz
number'' is used to describe the number of components of the Hurwitz space,
while what we call the Hurwitz number is called the ``degree''.
\end{rem}

\section{Reduction to four points}

The goal of this section is to use the machinery of limit linear series 
to prove:

\begin{prop}\label{reduce} To prove Theorem \ref{main} in general, it is 
enough to give a proof in the case that $r=4$.
\end{prop}

In order to use a degeneration argument for Proposition \ref{reduce}, the
key fact which we need (and which is lacking in the higher-genus case) is:

\begin{prop}\label{g0-dom} Every component of the space $G^1_d$ of 
Proposition \ref{equivs} maps dominantly under the forgetful map to 
$\cM_{0,r}$.
\end{prop}

\begin{proof} Indeed, we know \cite[Thm.\ 2.3]{e-h4} that if we fix 
ramification points, we have only finitely many $\fg^1_d$'s with the 
prescribed ramification, and that conversely, if we move the branch points, 
our rational function can always be deformed \cite[Cor.\ 3.2]{os3}; the
statement then follows by a dimension count, as in the proof of {\it ibid.}
\end{proof}

We will make essential use of the $r=1$ case of limit linear series,
developed by Eisenbud and Harris in \cite{e-h1}. We briefly review the
critical points of their theory in this case, where it becomes considerably
simpler. Suppose that we have a family $\cC$ of curves, with smooth generic
fiber, but with some nodal fibers. We assume that all nodal fibers are of
compact type, i.e., that their dual graph is a tree. Eisenbud and Harris 
construct a space over all of $\cC$ which correspond to usual $\fg^1_d$'s 
on smooth fibers of $\cC$, but correspond to {\it limit linear series} on 
the nodal fibers; by abuse of notation, we write $\fg^1_d$ to mean also
limit linear series. Suppose that $C$ is a nodal fiber
with (necessarily smooth) components $C_1,\dots,C_m$. In our case of $r=1$, 
a (refined) limit linear series on $C$ may be expressed as an $m$-tuple
of {\it aspects} on each $C_i$, where an aspect is a $\fg^1_{d_i}$ with
$d_i \leq d$, and the sole compatibility condition is that if $C_i$ and
$C_j$ meet at a node $P$, then the ramification index at $P$ of the 
aspects on $C_i$ and $C_j$ should be the same. Given $r$ smooth
sections $P_i$ of $\cC$, the Eisenbud-Harris construction also works to 
give spaces of $\fg^1_d$'s with at least a specified amount of ramification 
at the $P_i$ (in fact, limit linear series should in general allow for
base points away from the nodes, but since we will work with the case
that all ramification is specified, this won't arise).

We review the situation further in the case $g=0$, with all ramification
specified. This is studied in \cite[Thm.\ 2.4]{os2}; there, the families
considered involve only breaking off one component at a time, but our 
assertions here easily 
follow by the same arguments. For the rest of the section, we fix our
degenerate curve:

\begin{sit} The curve $C_0$ is the totally degenerate curve given by a
nodal chain of $r-2$ copies of $\PP^1$, with $P_1,P_2$ on the first 
component, $P_i$ on the $(i-1)$st component for $i<2<r-1$, and 
$P_{r-1},P_r$ on the last component.
\end{sit}

We consider families $\cC$ near a fiber isomorphic to the specified $C_0$.
Because all ramification is specified, the space of $\fg^1_d$'s
is finite over $\cC$, and is in fact finite \'etale in a neighborhood of 
$C_0$. Furthermore, a $\fg^1_d$ on $C_0$ is uniquely described by a 
collection of ramification indices $(e'_2,\dots,e'_{r-2})$ at the nodes,
which are required to satisfy a collection of triangle
inequalities and a parity condition. Specifically, if we consider any 
consecutive triple $e,e',e''$ starting with an odd-indexed term in the 
sequence 
$$e_1,e_2,e'_2,e_3,\dots,e_{r-2},e'_{r-2},e_{r-1},e_r,$$
we need to have $e \leq e'+e''$, $e' \leq e+e''$, and $e'' \leq e+e'$, and
we need $e+e'+e''$ to be odd. 

For later use, we note that the second condition implies immediately that
the triangle inequalities are in fact always strict, and also that
the allowed parity of $e'_2,\dots,e'_{r-2}$ is fixed by the $e_i$.

With these tools in hand, we can now complete our geometric argument.

\begin{proof}[Proof of Proposition \ref{reduce}] We fix the totally 
degenerate curve $C_0$ as in the above situation, and work with a
local universal family $\cC$ of genus-0 curves in a neighborhood of $C_0$,
denoting the generic curve of this family (which is also the generic
curve of $\cM_{0,r}$) by $C_{\eta}$.
It is enough to show that the relative $G^1_d$ space (with the desired
ramification at the marked points) is irreducible over the family
$\cC$, since by the previous proposition, every component of the global 
$G^1_d$ space meets the generic curve $C_{\eta}$. By the same token, it
is enough to show that any two $\fg^1_d$'s on the geometric generic fiber
$\bar{C}_{\eta}$ lie on the same
component of $G^1_d$. Furthermore, because the space of $\fg^1_d$'s is 
reduced over $C_0$, we have cannot have two components of $G^1_d$ meet over
$C_0$, so over our family $\cC$, irreducible components of $G^1_d$ are the
same as connected components. 

Accordingly, suppose we are given two $\fg^1_d$'s on $\bar{C}_{\eta}$. 
By the above discussion, these can be specialized to $\fg^1_d$'s
on $C_0$, which are described by the data of ramification indices
$(e'_2,\dots,e'_{r-2})$ and $(e''_2,\dots,e''_{r-2})$
respectively. We set the convention that $e'_1=e''_1:=e_1$,
and $e'_{r-1}=e''_{r-1}:=e_{r-1}$. Our claim is as follows: 
if we assume the $r=4$ case of
Theorem \ref{main}, then any two $\fg^1_d$'s on $C_0$ such that
$e'_i=e''_i$ for all but one $i$ necessarily lie on the same
component of $G^1_d$. 

Indeed, if we fix a node of $C_0$ corresponding to $e'_i$ (i.e., the 
$(i-1)$st node), we can restrict the family
$\cC$ to the closed subfamily $\cC_i$ in which only the chosen node of $C_0$
is allowed to be smoothed, giving a smooth component containing the two 
marked points $P_i$ and $P_{i+1}$, and the $(i-2)$nd and $i$th nodes (unless
$i=2$ or $r-2$, in which case $P_1$ or $P_r$ takes the place of the 
$(i-2)$nd or $i$th node respectively). The other components
remain fixed, so we may consider $\cC_i$ to be obtained from
the universal family over $\overline{\cM}_{0,4}$ by localizing around a 
degenerate curve, and gluing appropriate chains of $\PP^1$'s at the first 
and fourth marked points; in particular, the base of $\cC_i$ is naturally
a local scheme $U$ of $\overline{\cM}_{0,4}$ at a point corresponding to a
degenerate curve. If we write $\cC_{0,4}$ for the universal curve over
$U$, the point is to relate the $G^1_d$ spaces associated to $\cC_{0,4}$ 
and $\cC_i$.

Specifically, suppose we have chosen indices $e'_j=e''_j$ for all $j \neq i$.
For the sake of clarity, we denote by $G^1_d(\cC)$ our original space
of $\fg^1_d$'s on $\cC$, and by $G^1_d(\cC_i)$ and $G^1_d(\cC_{0,4})$ the 
spaces of $\fg^1_d$'s on $\cC_i$ and $\cC_{0,4}$. For the first two
spaces, we impose ramification $e_i$ at each $P_i$, so that $G^1_d(\cC_i)$
is simply the base change of $G^1_d(\cC)$, while for $G^1_d(\cC_{0,4})$
we impose ramification $e'_{i-1},e_{i},e_{i+1},e'_{i+1}$ at the four marked
points. Now, if we consider the closed subscheme $Z_i$ of
$G^1_d(\cC_i)$ which corresponds to limit $\fg^1_d$'s with
ramification indices $e'_j$ at the nodes (for $j \neq i$), the limit 
$\fg^1_d$'s are
uniquely determined except on the component with four marked points, so
$Z_i$ is isomorphic to the space $G^1_d(\cC_{0,4})$ which we have described.
Thus if we assume Theorem \ref{main} in the case $r=4$, we see
that the subscheme $Z_i$ of $G^1_d(\cC_i)$ is irreducible, so that 
any two $\fg^1_d$'s on $C_0$ for which $e'_i=e''_i$ for all but 
one $i$ lie on the same connected component of $G^1_d(\cC_i)$, 
and hence of $G^1_d(\cC)$. 

This proves the claim, and since every limit $\fg^1_d$ on $C_0$ can be 
smoothed to a $\fg^1_d$ on $\bar{C}_{\eta}$, the following numerical lemma 
completes the proof of our proposition.
\end{proof}

\begin{lem} Let $C_0$ be a totally degenerate marked curve of genus $0$,
and suppose we are given two $\fg^1_d$'s with ramification indices $e_i$
at the marked points, and classified by ramification indices 
$(e'_2,\dots,e'_{r-2})$ and $(e''_2,\dots,e''_{r-2})$
respectively at the nodes. Then it is possible to modify
$(e'_2,\dots,e'_{r-2})$ into $(e''_2,\dots,e''_{r-2})$,
by a sequence of changes affecting only one index at a time, and with 
every intermediate set of indices corresponding to a valid $\fg^1_d$ on 
$C_0$.
\end{lem}

\begin{proof} 
Suppose we have a $\fg^1_d$ on $C_0$ specified by the set 
$(e'_2,\dots,e'_{r-2})$. Since the allowed parity of each of
$e'_2,\dots,e'_{r-2}$ is fixed by the $e_i$, as 
long as we change them by $2$ at a time, we do not need to worry about 
violating the parity condition. It is thus enough to show that if 
$(e'_2,\dots,e'_{r-2})$ and $(e''_2,\dots,e''_{r-2})$
are distinct, there is always some $i$ with $e'_i \neq e''_i$
and for which we can increase $e'_i$ or $e''_i$ to make it closer 
to the other without violating any triangle inequalities. 
We prove this by induction.

We will induct on the following statement: suppose we are given $i$ such 
that 
$e''_i-e'_i\geq e''_{i-1} - e'_{i-1}$ and
$e'_i+2\leq e'_{i-1}+e_i$. Then either we can increase 
$e'_i$, or we must have
$e''_{i+1}-e'_{i+1}\geq e''_{i} - e'_{i}$ and 
$e'_{i+1}+2\leq e'_{i}+e_{i+1}$.
Indeed, if we cannot increase  
$e'_{i}$, the only triangle inequalities that could be violated are 
$e'_{i}+2 \leq e'_{i-1}+e_{i}$ or 
$e'_{i}+2 \leq e_{i+1}+e'_{i+1}$. But the first one is satisfied
by hypothesis, so the only possibility is that
$e'_{i}+2 > e_{i+1}+e'_{i+1}$, in which case we see we must have 
$e'_{i}+1=e_{i+1}+e'_{i+1}$. But we then see that
$$e''_{i+1}-e'_{i+1}=e''_{i+1}+e_{i+1}-e'_i-1 \geq
e''_i-e'_i$$
by the triangle inequality. Furthermore,
$e'_{i+1}+2 \leq e'_i+e_{i+1}$ because $e_{i+1} \geq 2$.

Now suppose that $i_0$ is the smallest number with 
$e'_{i_0} \neq e''_{i_0}$. Without loss of generality, we may
assume that $e'_{i_0}<e''_{i_0}$. But we see that this satisfies
the hypotheses of our inductive statement: the first inequality is clear
since $e''_{i_0-1}=e'_{i_0-1}$, while the second follows because we
have $e'_{i_0}+2 \leq e''_{i_0} \leq e''_{i_0-1}+e_{i_0}=
e'_{i_0-1}+e_{i_0}$. 
But by induction, we see that we must eventually be able to increment one 
of the $e'_i$ for $i \geq i_0$, since when $i=r-2$, we have 
$e''_{r-1}=e'_{r-1}=e_{r-1}$. This proves the lemma.
\end{proof}

\section{The case of four points}

In this section, we study the case of four points from the group-theoretic
point of view. Our setup throughout this section is as follows:

\begin{sit}\label{sit-4pts} We are given $d>0$, and 
$\vec{e}:=(e_1,e_2,e_3,e_4)$, with
$2d-2=\sum_i(e_i-1)$, and $2 \leq e_1 \leq e_2 \leq e_3 \leq e_4 \leq d$.
\end{sit}

We observe for later use that in our situation, we have 
$e_1+e_3 \leq d+1$, $e_2+e_4 \geq d+1$, $e_1+e_2 \leq d+1$, and 
$e_3+e_4 \geq d+1$. The first two inequalities follow from 
$e_1+e_3 \leq e_2+e_4$ together with $e_1+e_2+e_3+e_4=2d+2$, while the
second two follow by comparing with the first two.

Throughout this section, we will write sequences of the form $i,i+1,\dots,j$
(and similarly for descending sequences). If $j \geq i$, the meaning is clear:
an ascending sequence of length $j-i+1$. However, without further comment we
will also allow $j=i-1$, in which case the meaning will be the empty 
sequence (still of length $j-i+1$).

Our main result is the following:

\begin{thm}\label{4pts-thm} In Situation \ref{sit-4pts}, the Hurwitz number 
$h(d,4,0,\vec{e})$ is given by $\min\{e_i(d+1-e_i)\}_i$. 

Moreover, the possible Hurwitz factorizations
$(\sigma_1,\sigma_2,\sigma_3,\sigma_4)$ are classified explicitly as 
follows:
\begin{ilist}
\itm if $\sigma_3 \sigma_4$ is trivial or a single cycle, then we have 
\begin{multline*}\sigma_1=(d,d-1,\dots,e_3+e_4+1-k, \\
\sigma^{-(d+2-k-e_1)}(\ell),\sigma^{-(d+3-k-e_1)}(\ell),\dots,
\sigma^{-(e_3+e_4+1-2k)}(\ell)=\ell), \\
\sigma_2=(e_3+e_4+1-k,e_3+e_4+2-k,\dots,d-1,d,\ell,
\sigma^{-1}(\ell),\dots,\sigma^{-(d+1-k-e_1)}(\ell)), \\
\sigma_3=(k,k-1,\dots,2,1,e_4+1,e_4+2,\dots,e_3+e_4-k), \\
\sigma_4=(1,\dots,e_4),\end{multline*}
where we allow any $k$ with $e_3+e_4-d \leq k \leq e_3$
and $k \leq d+1-e_2$,
we define $\sigma:= (k,k+1,\dots,e_3+e_4-k)=\sigma_3 \sigma_4$,
and for a given $k$, we allow $\ell$ to vary in the range 
$k \leq \ell \leq e_3+e_4-k$.
\itm if $\sigma_3 \sigma_4$ is a product of two disjoint cycles, then we have
\begin{multline*}\sigma_1=(m+e_1-1,m+e_1-2,\dots,m+1,m), \\
\sigma_2=(d,d-1,\dots,m+e_1,m+d+k-e_3-e_4,m+d-1+k-e_3-e_4, \dots,k), \\
\sigma_3=(k,k-1,\dots,1,e_4+1,e_4+2,\dots,m+e_1-1, \\
m,m-1,\dots, m+d+1+k-e_3-e_4, m+e_1,m+e_1+1,\dots,d), \\
\sigma_4=(1,\dots,e_4),\end{multline*}
where we allow any $k$ with $1 \leq k \leq e_3+e_4-d-1$, and any $m$ with
$e_4-e_1+1 \leq m \leq d+1-e_1$ and $m \leq e_4$.
\end{ilist}
\end{thm}

Before giving the proof, we give a number of simple technical lemmas and 
their consequences; although each result individually is quite easy and
presumably well-known, we include them for the sake of staying as
self-contained as possible.

We begin by simplifying the transitivity condition on Hurwitz factorizations
in our situation.

\begin{lem}\label{repeat} Suppose that 
$\sigma_1,\sigma_2,\sigma_3,\sigma_4 \in S_d$ are cycles of length
$e_1,e_2,e_3,e_4$, with trivial product. Then the following are equivalent:
\begin{alist}
\itm the $\sigma_i$ form a Hurwitz factorization for $(d,4,0,\vec{e})$;
\itm every number in $\{1,\dots,d\}$ is in the support of at least
one of the $\sigma_i$; 
\itm every number in $\{1,\dots,d\}$ is in the support of exactly two
of the $\sigma_i$, except that either there is some $k$ which is in the 
support of every $\sigma_i$, or there exist $k \neq \ell$, with each 
in the support of three of the $\sigma_i$.
\end{alist}
\end{lem}

\begin{proof} It is clear that the transitivity condition for a) implies b). 
To see that b) implies c), the key point is that the identity 
$2d-2=\sum_i(e_i-1)$ leaves relatively little flexibility for the 
$\sigma_i$. Specifically,
if every number in $\{1,\dots,d\}$ is in the support of at least one cycle, 
it must be in the support of at least two, as otherwise the product could 
not be trivial. But we have $\sum_i e_i=2d+2$, and $2d$ of the numbers 
in the support of the $\sigma_i$ are accounted for, leaving only $2$ which 
could be in the support of more than two cycles. We conclude c).

Finally, to see that c) implies a), we need only check transitivity.
Since every number is in the support of at least two $\sigma_i$, we
cannot have any $\sigma_i$ disjoint from all the others. Thus, the only
way they could fail to generate a transitive subgroup would be if two of
the $\sigma_i$ were disjoint from the other two. But this cannot occur,
as the inequalities $e_3+e_4 \geq d+1$ and $e_2+e_4 \geq d+1$ imply that
$\sigma_4$ cannot be disjoint from either of $\sigma_2$ or $\sigma_3$.
\end{proof}

We next pursue a detailed study of the relationship between pairs of cycles
and their products.

\begin{lem}\label{two-cyc} Suppose
$\sigma_1, \sigma_2 \in S_d$ are non-disjoint cycles in
$S_d$, and let $\sigma$ be any cycle in the decomposition of
$\sigma_1 \sigma_2$ into disjoint cycles. Then there exists an element of
$\{1,\dots,d\}$ in the support of $\sigma$, $\sigma_1$, and $\sigma_2$.
\end{lem}

\begin{proof} This is routine: if
$\sigma$ consisted entirely of numbers in the support of $\sigma_1$ but not
$\sigma_2$, it would have to be equal to $\sigma_1$, contradicting the
non-disjointness hypothesis, and similarly with the $\sigma_i$ reversed.
One then verifies that to switch from elements in the support of $\sigma_1$
to elements in the support of $\sigma_2$ requires an element of $\sigma$
in the support of both. 
\end{proof}

\begin{lem}\label{cyc-prod-form} Let $\sigma,\sigma'$ be non-disjoint 
cycles, with $\sigma \sigma' \neq 1$.
Then there exists a unique expression (up to cycling of indices) of 
$\sigma'$ as $(w'_1,v'_1,\dots,w'_m,v'_m)$ and $\sigma$ as 
$(w_1,v_1,w_2,v_2,\dots,w_m,v_m)$ where the 
$w'_i,v'_i$ and $w_i,v_i$ are sequences of numbers, satisfying:
\begin{ilist}
\itm the $w_i$ and $w'_i$ are all non-empty, but the $v_i$ and $v'_i$ may
be empty;
\itm each $v'_i$ consists of numbers not in the support of $\sigma$; 
\itm each $v_i$ consists of numbers not in the support of $\sigma'$; 
\itm there exists a permutation $\tau \in S_m$ such that each $w_i$ is the 
inverse of $w_{\tau(i)}'$ (i.e., the same sequence in reversed order);
\itm if for all $i$ we set $k_i$ to be the first number in $w_i$, the set 
of $k_i$ is precisely the set of numbers in the support of all three of 
$\sigma,\sigma'$, and $\sigma \sigma'$.
\end{ilist}
\end{lem}

\begin{proof} By Lemma \ref{two-cyc}, there is some number in the support 
of $\sigma$, of $\sigma'$ and of $\sigma \sigma'$; we begin by designating
one such number to be $k_1$. In order to be able to write $\sigma$ in the
desired form, the order of the remaining $k_i$ are then uniquely determined.
Furthermore, we see that each $w_i$ must consist of the longest word in
$\sigma$ which starts with $k_i$, contains only numbers also in the support
of $\sigma'$, and does not contain $k_j$ for $j \neq i$. This uniquely
determines each $w_i$, and the $v_i$ are what remain. We can then do the
same for the $w'_i$ and $v'_i$, except that the $k_i$ could appear in a
different order in $\sigma'$, giving us the permutation $\tau$. It remains 
to check that these expressions have the desired properties, specifically
(ii), (iii), and (iv).

Note that if $n \neq k_i$ for any $i$ is any number in the support of 
$\sigma$ and $\sigma'$, since $n$ isn't in the support of $\sigma\sigma'$, 
then $\sigma'(n)=\sigma^{-1}(n)$, so we see that $\sigma'(n)$ must also 
be in the support of $\sigma$, immediately prior to $n$ in the cycle 
representation. Applying this 
inductively gives that all such $n$ appear in the $w_i$ in $\sigma$ and
in the $w'_i$ in $\sigma'$, and that each $w_{\tau(i)}'$ is inverse to $w_i$,
as desired.
\end{proof}

The following corollary is quite special to the case of at most two 
repetitions.

\begin{cor}\label{repeat-num-cycles} Let $\sigma,\sigma'$ be cycles, and
write $S \subseteq \{1,\dots,d\}$ for the intersection of the supports of
$\sigma$,$\sigma'$ and $\sigma \sigma'$. Suppose that either:
\begin{Ilist}
\itm $\sigma=\sigma_3$, $\sigma'=\sigma_4$, 
$(\sigma_1,\sigma_2,\sigma_3,\sigma_4)$
a Hurwitz factorization for $(d,4,0,\vec{e})$;
\itm $\sigma$ and $\sigma'$ are not disjoint, and $\#S \leq 2$.
\end{Ilist}
Then the number of disjoint cycles in $\sigma \sigma'$ is
equal to $\#S$ and is at most $2$, and there is exactly one element of $S$ 
in the support of each disjoint cycle.
\end{cor}

\begin{proof} We first observe that (I) implies (II). Indeed, $\sigma_3$
and $\sigma_4$ cannot be disjoint since $e_3+e_4 \geq d+1$, and we have
$\#S \leq 2$ by Lemma \ref{repeat} c), since any $k$ which occurs in the
support of $\sigma_3 \sigma_4$ must also occur in the support of $\sigma_1$
or $\sigma_2$ because of the trivial-product condition.

We next argue that (II) implies the conclusion of the corollary.
Lemma \ref{two-cyc} immediately handles the case $\#S \leq 1$.
For $\# S=2$, we apply the above lemma, noting first that in
this case we can always cycle indices so that $\tau=1$, i.e., each
$w'_i$ is the inverse of $w_i$. Then note that the formula 
$$(w_1,v_1,w_2,v_2)(w_1',v_1',w_2',v_2')=(k_1,v_1',v_2)(k_2,v_2',v_1)$$
holds regardless of whether the $v_i$ and $v_i'$ have non-zero length; 
since $k_1$ and $k_2$ are assumed to be in the support of
$\sigma \sigma'$, we see that it must consist of two disjoint cycles.
\end{proof}

We are now ready to give the proof of Theorem \ref{4pts-thm}. For the sake
of clarity, we break the proof into four propositions.

\begin{prop} Each of the possibilities enumerated in Theorem \ref{4pts-thm}
gives a valid Hurwitz factorization, and $\sigma_3 \sigma_4$ is in the
asserted form (and in particular, consists of the asserted number of cycles).
\end{prop}

\begin{proof} The main technicality is to check that the
individual cycles themselves make sense. This involves checking three
points: first, that all the elements listed lie in $\{1,\dots,d\}$;
second, that each word has non-negative length; and third, that there is
no overlap between the words in a given cycle. In fact, we first see that
all words have non-negative length, which then simplifies checking that 
their ranges are appropriate.

Before that, we observe that $\sigma=(k,k+1,\dots,e_3+e_4-k)$ makes sense:
the length is at least $1$, since we have $2k \leq 2e_3 \leq e_3+e_4$; 
and we check both $k \geq 1$ and $e_3+e_4-k \leq d$ using the inequality
$e_3+e_4-d \leq k$ (together with $e_3+e_4 \geq d+1$ for the first).

In general, we allow sequences to have length $0$, except those containing 
$k$, $\ell$, or $m$. 
In case (i), for $\sigma_1$ we require that $d \geq e_3+e_4-k$ and 
$e_3+e_4+1-2k \geq d+2-k-e_1$, which are equivalent to 
$k \geq e_3+e_4-d$ and $k \leq d+1-e_2$ respectively. For $\sigma_2$,
we need $e_3+e_4-k \leq d$ and $d+1-k-e_1 \geq 0$, giving 
$k \geq e_3+e_4-d$ and $k \leq d+1-e_1$ respectively. Since $d+1-e_2 \leq
d+1-e_1$, the last inequality will also be satisfied. Finally, for 
$\sigma_3$ we need $k \geq 1$ and $e_3+e_4-k \geq e_4$; the first is 
satisfied since $e_3+e_4-d \geq 1$, while the second is simply $k \leq e_3$.
Case (ii) is similar, with the only inequality appearing other than those
imposed directly being $m \geq e_3+e_4-d$. However, this is okay, 
since we have $m \geq e_4-e_1+1$, and the inequality $e_1+e_3 \leq d+1$ 
implies that $e_3+e_4-d \leq e_4-e_1+1$. Thus,
the ranges provided guarantee that the cycles make sense, and are in fact
equivalent to having the lengths of all words containing $k$, $\ell$, or $m$
be at least $1$, and the lengths of the remaining words being at least $0$.

We next address the first and third points simultaneously. In case (i),
everything follows easily from the ranges imposed for $k$: for $\sigma_4$ 
there is nothing to check; for $\sigma_3$ we check that $k < e_4+1$
and $e_3+e_4-k \leq d$; and for $\sigma_2$ and $\sigma_1$, everything is
immediate, since the terms involving $\sigma$ are automatically in the
correct range, and the sequence in $\sigma_1$ involving $\sigma$ could not 
wrap around without the sequence in $\sigma_2$ having negative length, and 
vice versa.

Similarly, in case (ii), the only points requiring any non-immediate checking 
are: for $\sigma_3$, that $k < m+d+1+k-e_3-e_4$, and $m<e_4+1$, with 
the former following from $e_1+e_3 \leq d+1$; for $\sigma_2$, that 
$m+d+k-e_3-e_4<m+e_1$; and for $\sigma_1$, that $m \geq 1$. Thus,
all the cycles consist of non-overlapping entries in $\{1,\dots,d\}$.
 
We can then check directly that the cycles are of the correct length and 
have trivial product, as well as that $\sigma=\sigma_3 \sigma_4$. Finally, 
using that b) implies a) in Lemma \ref{repeat} makes it easy to check 
that the cycles generate transitive subgroups of $S_d$, so all the 
possibilities listed are valid Hurwitz factorizations. 

It remains only to note that in case (i), we have already written $\sigma$
explicitly, so we see that $\sigma_3 \sigma_4$ is in fact trivial or a
single cycle, while in case (ii), we check that $m>m+d+k-e_3-e_4$, so that
$\sigma_1$ is disjoint from $\sigma_2$, and since 
$\sigma_1 \sigma_2 \sigma_3 \sigma_4=1$, it follows that $\sigma_3 \sigma_4$
is a product of two disjoint cycles.
\end{proof}

\begin{prop} No two possibilities enumerated in Theorem \ref{4pts-thm} are
equivalent.
\end{prop}

\begin{proof}
Cases (i) and (ii) of Theorem \ref{4pts-thm} are clearly invariant under 
relabeling. In case 
(i), we see that $k$ is determined as the number of elements in the support 
of both $\sigma_3$ and $\sigma_4$, so is invariant under relabeling. 
If $\sigma=1$ (i.e., if $k=e_3+e_4-k$), we have $\ell=k$ is the only 
possibility. Given $k$ with $\sigma \neq 1$, we see that $\ell$ is 
determined as the unique 
number (in the allowed range) such that $\sigma^{\ell-k}(k)$ is in the
support of $\sigma_1$ 
and $\sigma_2$ (and necessarily $\sigma$), so two possibilities with 
different $\ell$ cannot be equivalent.

In case (ii), the size of the intersection of the supports of $\sigma_1$ 
and $\sigma_4$ is $e_4+1-m$, so $m$ is relabeling-invariant. The overlap 
between the supports of $\sigma_3$ and $\sigma_4$ consists of two 
contiguous words, and $k$ is determined as the length of the word with
non-empty overlap with $\sigma_2$. Hence, no two possibilities are 
equivalent.
\end{proof}

\begin{prop} Every Hurwitz factorization is equivalent to one of the 
possibilities enumerated in Theorem \ref{4pts-thm}.
\end{prop}

\begin{proof} 
We begin by noting that by Corollary \ref{repeat-num-cycles}, we must
have that $\sigma_3\sigma_4$ consists of $0$, $1$, or $2$ disjoint cycles.
Furthermore, if $\sigma_3 \sigma_4=1$, then we
have $\sigma_3=\sigma_4^{-1}$, and $\sigma_1=\sigma_2^{-1}$, and 
$e_1=e_2=d+1-e_3=d+1-e_4$, so it is easy to check that the only possibility
is the $k=\ell=e_3=e_4$ case of (i). We can thus assume that 
$\sigma_3 \sigma_4 \neq 1$.

The first case we consider is that $\sigma=\sigma_3 \sigma_4$ is a single
cycle, or, equivalently by Corollary \ref{repeat-num-cycles}, that
there is a single number 
$k' \in \{1,\dots,d\}$ which is in the support of $\sigma_3,\sigma_4$, and 
in $\sigma=\sigma_3 \sigma_4$. 
Let $k$ be the number of elements in the support of both $\sigma_3$ and
$\sigma_4$. We may then relabel so that 
$\sigma_4=(1,\dots,e_4)$, and $k'$ gets mapped to $k$; i.e, so that 
the unique number in the support of $\sigma_3,\sigma_4$, and $\sigma$ is $k$.
Applying Lemma \ref{cyc-prod-form} to $\sigma_3$ and $\sigma_4$
with the only $k_i$ being
$k$ gives us that $\sigma_3$ is necessarily of the form 
$(k,k-1,\dots,2,1,a_1,\dots, a_{e_3-k})$ for some 
$a_i \in \{e_4+1,\dots,d\}$; relabeling the latter
range allows us to put $\sigma_3$ in the desired form. 

Next, note that by Lemma \ref{repeat}, there must be a unique number $\ell$
in the support of $\sigma_1$, of $\sigma_2$, and of $\sigma$. 
We then have also by Lemma \ref{repeat} that all the numbers 
$\{e_3+e_4-k+1,\dots,d\}$ must be in the support of $\sigma_2$, and 
we claim that they must be in a contiguous word, and followed immediately 
by $\ell,\sigma^{-1}(\ell),\dots,\sigma^{-(d+1-k-e_1)}(\ell)$. The claim 
is checked by applying Lemma \ref{cyc-prod-form} to $\sigma_2$ and 
$\sigma$, using that $\sigma_2 \sigma=\sigma_1^{-1}$, so that the only $k_i$
is $k_1=\ell$. The claim implies that we are free to reorder  
$\{e_3+e_4-k+1,\dots,d\}$ so that they appear in order, and furthermore so 
that $\sigma_2(d)=\ell$. Hence, we have put $\sigma_2$ in the desired 
form, and then $\sigma_1$ is determined by $\sigma_1 \sigma_2 \sigma=1$.

We next consider the case that $\sigma$ is a product of two disjoint cycles,
which by Corollary \ref{repeat-num-cycles} is equivalent to having 
two numbers $k',k'' \in \{1,\dots,d\}$ which are each in the support of 
$\sigma_3,\sigma_4$, and in $\sigma:=\sigma_3 \sigma_4$. Then $k'$ is in 
one of the disjoint cycles of $\sigma$, and $k''$ is in the other. By
Lemma \ref{repeat}, we see that since we already have $k',k''$ occurring in
$\sigma_3,\sigma_4$ and $\sigma$ (hence in either $\sigma_1$ or $\sigma_2$),
we cannot have any numbers occurring in $\sigma_1,\sigma_2$ and $\sigma$.
By Corollary \ref{repeat-num-cycles} (II), we see that $\sigma_1$ and 
$\sigma_2$ must be disjoint, and since $\sigma_1\sigma_2=\sigma^{-1}$, 
we see that $k'$ is in the 
support of one, and $k''$ is in the support of the other;
without loss of generality, we may assume that $k'$ is in the support of 
$\sigma_1$ and $k''$ in $\sigma_2$. We also note that this implies that
each of $\{1,\dots,d\}$ is in the support of either $\sigma_3$ or 
$\sigma_4$.

We once again normalize so that $\sigma_4=(1,\dots,e_4)$, and we can
further require that if we write $\sigma_3=(w_1,v_1,w_2,v_2)$
and $\sigma_4=(w_1',v_1',w_2',v_2')$ as in Lemma \ref{cyc-prod-form}, we
can set $w_1'=(1,2,\dots,k)$, with $k$ being the corresponding relabeling
of $k''$, i.e., the unique number in the support of $\sigma_3$,$\sigma_4$,
and $\sigma_2$. We then have $w_1$ in the desired form, and $w_2$ will
likewise be in the desired form for some $m$, which will necessarily be 
the unique number in the support of $\sigma_3$,$\sigma_4$, and $\sigma_1$. 
Relabelling $e_4+1,\dots,d$ as necessary, we can place $v_1$ and $v_2$,
hence $\sigma_3$ in the desired form, and $\sigma_1$ and $\sigma_2$ are
then uniquely determined as disjoint cycles with 
$\sigma_1 \sigma_2 \sigma =1$, and containing $m$ and $k$ respectively.

This then completes the proof of the claim that every Hurwitz factorization
is equivalent to one of the enumerated possibilities.
\end{proof}

\begin{prop} The number of possibilities enumerated in Theorem \ref{4pts-thm}
is equal to $\min\{e_i(d+1-e_i)\}_i$.
\end{prop}

\begin{proof} The formula $\min\{e_i(d+1-e_i)\}_i$ falls into two 
situations: if $e_4 \geq d+1-e_1$, then it is equal to $e_4(d+1-e_4)$, 
while if $e_4 \leq d+1-e_1$, then it gives $e_1(d+1-e_1)$. 

We first consider the situation that $e_4 \geq d+1-e_1$.
Here, because $e_4+e_1 \geq d+1$, we have $e_2+e_3 \leq d+1$,
so $e_3 \leq d+1-e_2$, and in case (i) of Theorem \ref{4pts-thm}
the inequality $e_3+e_4-d \leq k \leq e_3$ automatically implies
$k \leq d+1-e_2$. We thus have 
$$\sum_{k=e_3+e_4-d}^{e_3} \sum_{\ell=k}^{e_3+e_4-k} 1 =
\sum_{k=e_3+e_4-d}^{e_3} (e_3+e_4-2k+1) = (d+1-e_3)(d+1-e_4)$$
possibilities from case (i). Similarly, we have $d+1-e_1 \leq e_4$ 
so $e_4-e_1+1 \leq m \leq d+1-e_1$ implies that $m \leq e_4$.
Thus, our ranges are $1 \leq k \leq e_3+e_4-d-1$ and 
$e_4-e_1+1 \leq m \leq d+1-e_1$, yielding $(e_3+e_4-d-1)(d+1-e_4)$
possibilities in case (ii), and giving us the desired $e_4(d+1-e_4)$ 
possibilities in total (note that $e_3+e_4-d-1$ and $d+1-e_4$ are
always non-negative, so these formulas are always valid). 

The situation that $e_4 \leq d+1-e_1$ proceeds similarly, with $e_1e_2$ 
possibilities arising from case (i), and $e_1(d+1-e_1-e_2)$ possibilities
arising from case (ii).
\end{proof}

Combining the statements of the four propositions, we immediately conclude
Theorem \ref{4pts-thm}.

From the theorem, we deduce quite directly:

\begin{cor}\label{4pts-cor} In Situation \ref{sit-4pts}, the Hurwitz 
factorizations for $(d,4,0,\vec{e})$ all lie in a single pure braid orbit.
\end{cor} 

\begin{proof} We first see that all the factorizations in case (i) of the
theorem are in a single pure braid orbit, and then show that any 
factorization in case (ii) is in the same braid orbit as some factorization 
in case (i). 

Suppose we start with $(\sigma_1,\sigma_2,\sigma_3,\sigma_4)$
corresponding to a given $k,\ell$ of case (i). Our first claim is that
if we replace $(\sigma_1,\sigma_2)$ by 
$(\sigma_2^{-1} \sigma_1 \sigma_2, 
\sigma_2^{-1} \sigma_1^{-1} \sigma_2 \sigma_1 \sigma_2=
\sigma \sigma_2 \sigma^{-1})$, we stay in case (i), leaving $k$ fixed, 
while replacing $\ell$ by $\sigma(\ell)$. The first part is clear, while
the assertion on $\ell$ is checked by direct computation, using that
since $\sigma_3,\sigma_4$ remain fixed, it is enough to see what happens
to $\sigma_2$. Thus, for a
given $k$, every possible $\ell$ is in the same braid orbit. 

To analyze the Hurwitz factorizations for different $k$, for each
$k$ we set $\ell=k$, where we have 
$\sigma_1=(d,d-1,\dots,e_3+e_4+1-k,d+1-e_2,d-e_2,\dots,k)$
and hence
\begin{multline*}\sigma_2 \sigma_3=\sigma_1^{-1}\sigma_4^{-1}=
(k,k-1,\dots,1,e_4,e_4-1,\dots,d+2-e_2,\\
e_3+e_4+1-k,e_3+e_4+2-k,\dots,d).\end{multline*}
We check that if we replace
$(\sigma_2,\sigma_3)$ by $(\sigma'_2,\sigma'_3):=
(\sigma_3^{-1} \sigma_2 \sigma_3, 
\sigma_3^{-1} \sigma_2^{-1} \sigma_3 \sigma_2 \sigma_3)$, then as long
as $k$ is not minimal, we 
remain in case (i), but replace $k$ by $k-1$. Here, a relabeling is in
principle necessary, but we can instead check that $\sigma'_3 \sigma_4$ is 
still a single cycle, so that we remain in case (i), and that
the supports of $\sigma'_3$ and $\sigma_4$ overlap in $k-1$ elements.
We therefore see that every possibility
in (i) is always in a single pure braid orbit.

Finally, we suppose we have $(\sigma_1,\sigma_2,\sigma_3,\sigma_4)$
corresponding to a given $k,m$ of case (ii). In this case, we again
replace $(\sigma_2,\sigma_3)$ by $(\sigma'_2=\sigma_3^{-1} \sigma_2 \sigma_3,
\sigma'_3=\sigma_3^{-1} \sigma_2^{-1} \sigma_3 \sigma_2 \sigma_3)$, 
and note that since $\sigma_1,\sigma_4$ remain unchanged, $\sigma'_2$
determines $\sigma'_3$. One then computes that as long as $k<e_3+e_4-d-1$,
$\sigma'_2$ is still a possibility from case (ii), with $m$ the same,
but $k+1$ instead of $k$. Finally, if $k=e_3+e_4-d-1$, one checks that
applying the same pure braid operation, we move into case (i), with
$k=e_3+e_4+d$ (and $\ell=m$). Thus, every possibility in case (ii) is in 
the same pure braid orbit as some possibility in case (i), and we get that 
everything is in the same pure braid orbit.
\end{proof}

Using Proposition \ref{equivs}, and Proposition \ref{reduce}, we see 
immediately that Corollary \ref{4pts-cor} implies Theorem \ref{main},
and we are done.

\section{Loose ends}

We begin with a further remark in the case of four points. The Hurwitz
number $\min\{e_i(d+1-e_i)\}_i$ computes the number of rational functions
$\PP^1 \to \PP^1$ with four fixed branch points on the target. If instead 
we look at fixed ramification points on the source, we find that the number
is $\min\{e_i,d+1-e_i\}_i$ \cite[Rem.\ 5.9]{os2}. Despite the close 
geometric relationship between these two
numbers, there is no {\it a priori} reason for there to be any numerical
relationship at all, so their similarity is striking. We note further that
with the exception of the case that we have $e_i=d$ for some $d$, both
formulas are symmetric with respect to replacing the $e_i$ by $d+1-e_i$;
this motivates us to ask:

\begin{ques} Is there a natural involution on the set of rational
functions of degree $d$ having exactly four ramification points, which
replaces the ramification indices $e_i$ by $d+1-e_i$, and holds both the 
ramification and branch points fixed?
\end{ques} 

A more obvious question left unanswered by our analysis is:

\begin{ques} Is there a closed form for the genus-$0$ pure-cycle Hurwitz 
numbers for any number of branch points?
\end{ques}

Next, we observe that it is a consequence of Theorem \ref{main} that if
we fix $d,r$ and $\vec{e}$, all possible Hurwitz factorizations are in a
single Nielsen class, i.e., they generate
the same group, and are in the same conjugacy classes within that group.
However, with some non-trivial group theory and sufficient perseverance,
one can already see this quite directly:

\begin{thm}\label{nielsen} Given $d,r$ and $\vec{e}=(e_1,\dots,e_r)$ with 
$2d-2=\sum_i(e_i-1)$, and all $e_i\geq 2$, suppose we have 
$(\sigma_1,\dots,\sigma_r)$ and
$(\sigma'_1,\dots,\sigma'_r)$ two Hurwitz factorizations for
$(d,r,0,\vec{e})$, generating 
groups $G,G' \subseteq S_d$. Then there exists a simultaneous conjugation in 
$S_d$ making $G'=G$, and each $\sigma_i$ conjugate to $\sigma'_i$ inside $G$.
That is, any two Hurwitz factorizations lie in the same Nielsen class.

In fact, if $r=2$, we have $G$ isomorphic to the cyclic group $C_d$.
If $r=3$ with $(e_1,e_2,e_3)=(4,4,5)$, we have $G \cong S_5$, imbedded
as a doubly transitive subgroup of $S_6$. Otherwise, we always have 
$G=S_d$ or $G=A_d$ depending on the parity of the $e_i$. 
\end{thm}
 
\begin{proof} The case that $r=2$ is clear, as we must have $e_1=e_2=d$.

For $r=3$, we note that the first assertion is clear, since the Hurwitz
number is equal to $1$ by Lemma \ref{3pts}.

For $r > 3$, we reduce the first assertion to the second. In the case
that $G=S_d$, this is trivial, while in the case that $G=A_d$, we need only
observe that since $2d-2=\sum_i (e_i-1)$, and all $e_i \geq 3$, we can have 
at most one cycle of order greater than $d-2$.
We can always fix this cycle by simultaneous conjugation in $S_d$, and then
any cycles of given length less than or equal to $d-2$ are in the same 
conjugacy class in $A_d$.

For the second assertion, we begin by arguing that with $r>2$, we must have
$G$ primitive, 
i.e., that there is no non-trivial partition of $\{1,\dots,d\}$ into 
blocks on which the action of $G$ is well-defined.
Indeed, if there were such a partition, since $G$ is transitive the
blocks would all have to have the same size $m$, for some $m|d$.
We would then necessarily have each $\sigma_i$ either of size a
multiple of $m$, acting as a $e_i':=\frac{e_i}{m}$-cycle $\sigma'_i$ on 
$d':=\frac{d}{m}$ blocks of size $m$, or of size strictly less than
$m$, acting trivially on the blocks. Say we have $s$ of the latter;
without loss of generality, we may assume that $e_1,\dots,e_s<m$,
and $e_{s+1},\dots,e_r \geq m$. Then $\sigma_{s+1}',\dots,\sigma'_r$
give a Hurwitz factorization in $S_{d'}$, so we must have
$2d'-2 \leq \sum_{i=s+1}^r (e'_i-1)$. On the other hand, we compute
that since $2d-2=\sum_i (e_i-1)$, we have $2d+r-2=\sum_i e_i$, so
$$2\frac{d}{m}+\frac{r-2}{m}-\sum_{i=1}^s \frac{e_i}{m}=
\sum_{i=s+1}^r \frac{e_i}{m}=\sum_{i=s+1}^r e'_i,$$ 
and so $2d'-2-\delta=\sum_{i=s+1}^r (e'_i-1)$, where
$$\delta=\sum_{i=1}^s \frac{e_i}{m}-s+r-2-\frac{r-2}{m} \geq
\frac{2s}{m}-s+r-2-\frac{r-2}{m}=\frac{(m-1)(r-s-2)}{m}+\frac{s}{m},$$
so we must have $\frac{(m-1)(r-s-2)}{m}+\frac{s}{m} \leq 0$.
Since the $\sigma'_i$ act transitively on $d'$ elements, and have
trivial product, there must be at least $2$ of them which are non-trivial, 
so that $r-s-2 \geq 0$. Since $m>1$, we see that $\delta \geq 0$, and we can
have $\delta=0$ only if $r-s-2=s=0$, i.e., $r=2$. Thus, with our hypothesis
that $r>2$, we must have $\delta>0$, a contradiction. 

We note that in the case that $d \leq 3$, the only transitive subgroups are
$A_d$ and $S_d$, so there is nothing to prove. In the case $d=4$, one checks
directly that there is no primitive subgroup other than $S_4$
and $A_4$, so we need only consider the case $d \geq 5$.

Now, we wish to apply the theorem of Williamson \cite{wi2} stating that if 
a primitive subgroup of $S_d$ contains a cycle of order $e$, 
with $e\leq (d-e)!$, then it must be either $A_d$ or $S_d$.
Since we have $2d-2=\sum_i (e_i-1)$, we see that we must have
$e_i \leq \lfloor{\frac{2d-2}{r}+1\rfloor}$ for some $i$. One then computes 
directly that Williamson's theorem gives the desired result unless we have
$r=3, d \leq 10$, or $r=4,d \leq 5$. More specifically, the only cases
falling outside Williamson's theorem are $r=3$ with $(e_1,e_2,e_3)=
(3,4,4),(4,4,5),(5,5,5),(7,7,7)$ or $r=4$ with $(e_1,e_2,e_3,e_4)
=(3,3,3,3)$. In these cases, one can check directly that the group is
$A_d$ or $S_d$, as appropriate, except in the $(4,4,5)$ case, where one
can compute the group explicitly, checking that it is doubly transitive
and has order $120$, which is well-known to determine it uniquely.
\end{proof}

Our result is sharp in the sense that if one drops either the pure-cycle or
the genus-$0$ hypothesis, there are many examples for which the Hurwitz
space is not irreducible. However, there are nonetheless many examples
for which the Hurwitz space is irreducible which are not covered by our
main theorem. We will consider here one generalization which 
remains in the pure-cycle case, but seeks to drop the genus-$0$ hypothesis
in favor of an assumption that could be viewed philosophically as an 
effective form of the results of Conway-Fried-Parker-V\"olklein, in that 
it requires at least
$3g$ transpositions in order to apply. However, our result will be 
conditional on a positive answer to a geometric question, which we now 
discuss.

Zariski asked whether every Hurwitz space of genus-$g$ covers of $\PP^1$ 
with prescribed branching type over at least $3g$ points maps dominantly
to $\cM_g$ under the forgetful map. This is now known to be false in some
cases, but we will be interested in an analogous yet different 
question which arises when one wants to compare the points of view of 
linear series and branched covers:

\begin{ques}\label{zar-q} Fix $r,g \geq 0$, $d \geq 1$ and 
$\vec{e}=(e_1,\dots,e_r)$ with $2 \leq e_i \leq d$ for all $i$, and 
$2d-2-g=\sum_i (e_i-1)$. Consider the space 
$MR$ parametrizing tuples consisting of a genus-$g$ curve $C$, 
points $P_1,\dots, P_{r+3g}$ on $C$, and 
a map $f:C \to \PP^1$ of degree $d$, ramified to order $e_i$ at $P_i$ for 
$i \leq r$ and simply ramified at $P_{r+1},\dots,P_{r+3g}$. Does every 
component of $MR$ map dominantly to $\cM_{g,r}$ under the map induced by 
forgetting $f$ and $P_{r+1},\dots,P_{r+3g}$?
\end{ques}

The positive answer to this question in the case $g=0$ is Proposition 
\ref{g0-dom}. We also remark that Steffen \cite{st3} (see also \cite{ge1})
has a result along 
these lines for linear series of any degree and dimension, but without 
any ramification specified. He accomplishes this by studying degeneracy
loci of suitable maps of vector bundles; one might try to study our
question by looking at Schubert conditions on maps of vector bundles,
and suitable intersections of such conditions.

The application of Question \ref{zar-q} to irreducibility of Hurwitz spaces
is as follows.

\begin{thm}\label{high-gen} Fix $r,g,d,$ and $\vec{e}$ as above. Then a 
positive answer to Question \ref{zar-q} implies that $\cH(d,r,g,\vec{e})$
is irreducible, where $\cH(d,r,g,\vec{e})$ is the Hurwitz space of 
covers of $\PP^1$ of genus $g$ and degree $d$, with a single ramified point 
of index $e_i$ over the $i$th branch point for $i \leq r$, and simple 
branching over the remaining branch points. Equivalently, the set of
Hurwitz factorizations consisting of $e_i$-cycles and $3g$ transpositions
all lie in a single pure braid orbit.
\end{thm}

\begin{proof} We first consider the generalization of Proposition 
\ref{equivs} in this case. The argument for the equivalence of (i) and (ii)
goes through unmodified in the generality of higher-genus covers. The 
argument for the equivalence of (iii) and (iv), where in both cases we 
prescribe simple ramification at $3g$ additional unspecified points, is
likewise the same as in the genus $0$ case. We then have that the Hurwitz
space is the image of (a dense open subset of) $MR$, so we see that (iii) 
or (iv) imply (i) and (ii),
and it is enough to check (iv), i.e., to work from the point of view of
linear series.

A positive answer to Question \ref{zar-q} takes the place of Proposition
\ref{g0-dom}, and allows us to work over the generic $r$-marked curve of
genus $g$, or more specifically, locally around a given degenerate curve,
as in the genus $0$ case. Instead of working with a totally degenerate 
curve, we work with a curve $C_0$ consisting of a copy of $\PP^1$ with 
$r$ marked points, and with $g$ elliptic tails. As in the proof of 
\cite[Thm.\ 2.6]{os2}, the limit linear series on this curve are completely
determined by their aspects on $\PP^1$; on each elliptic tail, they 
consist of the degree $2$ map to $\PP^1$, simply ramified at the node (and
at three other points, which are uniquely determined as differing from the
node by $2$-torsion points). Furthermore, the ramification imposed at each
node on $\PP^1$ is simple ramification; thus, the limit linear series are
in natural bijection with the linear series on $\PP^1$ with the prescribed
ramification at $r+g$ points. We know by Theorem \ref{main} that the space
of these linear series is irreducible as we allow the $r+g$ ramification
points to move, so we conclude irreducibility of the space of $\fg^1_d$'s
in a neighborhood of $C_0$, and in particular, on the generic $r$-marked
curve of genus $g$, as desired.
\end{proof}

Results of Conway-Fried-Parker-V\"olklein (see \cite[Appendix]{f-v1}, and also
\cite{fr4}) show that, roughly speaking,
for any given group and collection of conjugacy classes, if every 
conjugacy class is repeated often enough, then the components of the 
Hurwitz space are determined by a certain invariant, called the lifting 
invariant. Our results fit into the same general philosophy, and might be 
thought of as an effective version of Conway-Fried-Parker-V\"olklein for the 
pure-cycle case.

Finally, we remark that Question \ref{zar-q} would potentially have 
interesting applications to the study of covers in positive characteristic, 
as well. One cannot hope for a positive answer outside characteristic $0$
without some further hypotheses: for instance, in the case $g=0$, the 
statement is known to fail if one does not require all $e_i<p$ (see
\cite[Ex.\ 5.6]{os7}). However, 
a positive answer in the case all $e_i<p$ would give an important step 
towards giving new non-existence results for tame covers in positive 
characteristic, as is carried out in the genus-$0$ case in \cite{os12}.

\bibliographystyle{hamsplain}
\bibliography{hgen}

\newcommand{\noopsort}[1]{} \newcommand{\printfirst}[2]{#1}
  \newcommand{\singleletter}[1]{#1} \newcommand{\switchargs}[2]{#2#1}
\providecommand{\bysame}{\leavevmode\hbox to3em{\hrulefill}\thinspace}
\begin{thebibliography}{10}

\bibitem{e-h4}
David Eisenbud and Joe Harris, \emph{Divisors on general curves and cuspidal
  rational curves}, Inventiones Mathematicae \textbf{74} (1983), 371--418.

\bibitem{e-h1}
\bysame, \emph{Limit linear series: Basic theory}, Inventiones Mathematicae
  \textbf{85} (1986), 337--371.

\bibitem{e-h7}
\bysame, \emph{Irreducibility and monodromy of some families of linear series},
  Annales scientifiques de l'\'Ecole Normale Sup\'erieure \textbf{20} (1987),
  no.~1, 65--87.

\bibitem{fr1}
Michael~D. Fried, \emph{Alternating groups and moduli space lifting
  invariants}, preprint.

\bibitem{fr4}
\bysame, \emph{Connectedness of families of sphere coverings of a given type},
  preprint.

\bibitem{fr3}
\bysame, \emph{Fields of definition of function fields and {H}urwitz families
  -- groups as {G}alois groups}, Communications in Algebra \textbf{5} (1977),
  no.~1, 17--82.

\bibitem{fr2}
\bysame, \emph{Relating two genus 0 problems of {J}ohn {T}hompson}, Progress in
  Galois theory, Developments in Mathematics, no.~12, 2005, pp.~51--85.

\bibitem{fr5}
\bysame, \emph{The {M}ain {C}onjecture of modular towers and its higher rank
  generalization}, Groupes de {G}alois arithmetiques et differentiels (Luminy
  2004), Seminaires et congres, vol.~13, 2006.

\bibitem{f-v1}
Michael~D. Fried and Helmut V\"olklein, \emph{The inverse {G}alois problem and
  rational points on moduli spaces}, Mathematische Annalen \textbf{290} (1991),
  no.~4, 771--800.

\bibitem{f-v2}
\bysame, \emph{The embedding problem over a {H}ilbertian {PAC}-field}, Annals
  of Mathematics \textbf{135} (1992), no.~2, 469--481.

\bibitem{fu3}
William Fulton, \emph{Hurwitz schemes and the irreducibility of moduli of
  algebraic curves}, Annals of Mathematics \textbf{90} (1969), 542--575.

\bibitem{ge1}
Joergen~Anders Geertsen, \emph{Push-forward of degeneracy classes and
  ampleness}, Proceedings of the AMS \textbf{129} (2001), no.~7, 1885--1890.

\bibitem{os12}
Brian Osserman, \emph{Linear series and existence of branched covers},
  \mbox{arXiv:math.AG/0507096}.

\bibitem{os2}
\bysame, \emph{The number of linear series on curves with given ramification},
  International Mathematics Research Notices \textbf{2003} (2003), no.~47,
  2513--2527.

\bibitem{os3}
\bysame, \emph{Deformations of covers, {B}rill-{N}oether theory, and wild
  ramification}, Mathematical Research Letters \textbf{12} (2005), no.~4,
  483--491.

\bibitem{os7}
\bysame, \emph{Rational functions with given ramification in characteristic
  $p$}, Compositio Mathematica \textbf{142} (2006), no.~2, 433--450,
  \mbox{arXiv:math.AG/0407445}.

\bibitem{st3}
Frauke Steffen, \emph{A generalized principal ideal theorem with an application
  to {B}rill-{N}oether theory}, Inventiones Mathematicae \textbf{132} (1998),
  no.~1, 73--89.

\bibitem{vo1}
Helmut V\"olklein, \emph{Groups as {G}alois groups}, Cambridge Studies in
  Advanced Mathematics, no.~53, Cambridge University Press, 1996.

\bibitem{wi2}
Alan Williamson, \emph{On primitive permutation groups containing a cycle},
  Mathematische Zeitschrift \textbf{130} (1973), 159--162.

\end{thebibliography}
\end{document}